\documentclass[conf]{new-aiaa} 
\usepackage[utf8]{inputenc}
\usepackage{textcomp}
\usepackage{graphicx}
\usepackage{amsmath}
\usepackage[version=4]{mhchem}
\usepackage{siunitx}
\usepackage{longtable,tabularx}
\usepackage{placeins}
\usepackage{booktabs} 
\usepackage{siunitx} 
\usepackage{float}
\usepackage{comment}
\usepackage{makecell}
\usepackage{multirow}
\usepackage{array}
\usepackage{bm}
\usepackage{tablefootnote}
\usepackage{threeparttable}
 \usepackage{gensymb}

\usepackage{graphicx}
\graphicspath{{images/}} 

\usepackage[nolist,nohyperlinks]{acronym}

\title{Logistics Analysis for Lunar Post-Mission Disposal}

\author{Evangelia Gkaravela\footnote{MSc Student, Department of Systems and Enterprises, AIAA Student Member.}}
\author{Hao Chen\footnote{Assistant Professor, Department of Systems and Enterprises, AIAA Member.}}
\affil{Stevens Institute of Technology, Hoboken, NJ, 07030}

\begin{document}

\maketitle

\begin{acronym}
    \acro{DSM}{Design Structure Matrix}
    \acro{DWE}{Direct Water Electrolysis}
    \acro{FSPS}{Fission Surface Power System} 
    \acro{HB}{Habitat}
    \acro{HR}{Hydrogen Reduction}
    \acro{ILRS}{International Lunar Research Station}
    \acro{IPEx}{ISRU Pilot Excavator}
    \acro{ISRU}{In-Situ Resource Utilization}
    \acro{LEO}{Low Earth Orbit}
    \acro{LCA}{Lifecycle Assessment}
    \acro{LCI}{Lifecycle Inventory}
    \acro{LCROSS}{Lunar Crater Observation and Sensing Satellite}
    \acro{LS}{Lunar Surface}
    \acro{MRE}{Molten Regolith Electrolysis}
    \acro{MCR}{Methane Carbothermal Reduction}
    \acro{NASA}{National Aeronautics and Space Administration}  
    \acro{SWE}{Soil/Water Extraction}  
\end{acronym}

\begin{abstract}
\lettrine{A}{s} human activities on the Moon expand through initiatives like NASA's Artemis program, the need for sustainable post-mission disposal strategies becomes critical to maintaining the lunar environment. This paper analyzes the logistics and environmental implications of waste products generated by \acl{ISRU} technologies employed in oxygen production on the Moon. The study examines the inputs, generation of products, and the resulting byproducts from \acl{MRE}, \acl{SWE}, and \acl{DWE} systems. These technologies yield varied byproducts, including slag, metals and volatiles, each presenting unique challenges for disposal and recycling. The analysis assesses the economic and ecological impacts of \acl{ISRU} activities on lunar operations using a multi-commodity flow model adapted from cislunar logistics frameworks. The results inform that ISRU-enabled missions achieve a significant threefold cost reduction. However, the management of byproducts remains a critical challenge, demanding innovative solutions to address their impact and support scalable and sustainable lunar exploration.


\end{abstract}

\section*{Nomenclature}

{\renewcommand\arraystretch{1.0}
\noindent\begin{longtable*}{@{}l @{\quad=\quad} l@{}}
$\mathcal{A}$  & set of directed arcs \\
$\mathcal{C}$ & set of commodities \\
$\mathcal{C_C}$ & set of continuous commodities \\
$\mathcal{C_D}$ & set of discrete commodities \\
$\bm{c}$   & cost coefficient vector \\
$\bm{d}$ & demand and supply vector \\
$\mathcal{G}$ & network graph \\
$g$  & binary interval variable \\
$Q$  & concurrency matrix \\
$i$  & node index ($\in \mathcal{N}$) \\
$\mathcal{J}$  & cost \\
$j$  & node index ($\in \mathcal{N}$) \\
$l$  & number of concurrency constraints \\
$\mathcal{N}$ & set of nodes \\
$Q$  & transformation matrix \\
$q$ & design quantity \\
$\mathcal{T}$ & set of time steps \\
$t$ & time index (integer) \\
$W$ & set of time windows \\
$\Delta t$ & time of flight \\
\end{longtable*}}

\section{Introduction}

\lettrine{T}{}he revival of lunar exploration, led by projects such as NASA’s Artemis, China’s \ac{ILRS}, and the European Space Agency’s Lunar Exploration missions \cite{zhang2023}, is a turning point in human endeavor for a sustained presence beyond Earth. NASA emphasizes that a sustained lunar presence will inspire humanity to explore even more distant worlds while fostering technological and scientific advancements that benefit life on Earth. \cite{nasa2020}.

The Moon's exploration has been a fundamental aspect of human space exploration, beginning with the Apollo and Luna missions, which offered initial insights into its geological features and surface conditions. These missions showcased human technological achievements and underscored the Moon's potential as a valuable resource for sustainable exploration \cite{anand2010, badescu2012}. Its proximity to Earth makes it an ideal platform for testing advanced technologies, such as habitat construction and life-support systems, which are essential for interplanetary missions. As a stepping stone to Mars and beyond, the Moon offers unparalleled opportunities for developing sustainable exploration frameworks.

While seemingly barren, the Moon holds unique scientific value as a natural laboratory for studying planetary evolution and serves as a platform for observing deep-space phenomena \cite{badescu2012}. Preserving its environment is critical not only for maintaining its scientific integrity but also for ensuring that future missions remain viable. Mainly, Lunar surface sustainability involves systems and infrastructure to support long-term human and robotic operations that do not degrade the Lunar environment or interfere with future missions.\cite{nasa2020}.

Current space missions rely on open-loop life support systems, which require transporting supplies from Earth and discarding waste. However, future lunar missions will depend on closed-loop systems to recycle waste streams, reducing reliance on Earth-based inputs and minimizing waste generation \cite{hightower1992}. When integrated with \ac{ISRU} technologies, these systems offer a pathway toward more sustainable lunar operations.
Particularly, \ac{ISRU} integration into mission architectures enables the production of mission-critical consumables, including propulsion fuels, power, and life support resources, that significantly reduce the mass, cost, and risk associated with missions.\cite{sanders2012,anand2010} 

\ac{ISRU} provides a promising approach to sustainable lunar exploration by facilitating the extraction and use of abundant local resources, such as regolith, which contains approximately 44\% oxygen along with oxides, and other volatiles, as well as water ice located in permanently shadowed craters \cite{badescu2012, basilevsky2012}. \ac{ISRU} techniques, including \acf{MRE}, \acf{SWE}, and \acf{DWE}, reduce reliance on Earth-based resupply by producing essential consumables such as oxygen, water, and hydrogen. 

Further postreactor processing of byproducts, including molten iron, silicon, aluminum, titanium, and glassy slag, enables the production of infrastructure \cite{badescu2012}, spare parts, and even solar arrays directly on the lunar surface, significantly enhancing the sustainability of operations \cite{landis2005, curreri2006}. However, these processes generate waste materials that, if not managed effectively, could harm the lunar environment. Issues such as dust generation, landscape alterations, and lunar surface contamination highlight the need for systematic waste disposal strategies. \cite{dallas2021}.

This study addresses the environmental and logistical challenges associated with \ac{ISRU} activities by analyzing the waste profiles of three key ISRU methods: \ac{MRE},  \ac{SWE}, and \ac{DWE}. While ISRU significantly reduces mission costs, its environmental impacts must be carefully quantified and weighed against potential cost savings. To achieve this, the study employs a multicommodity flow model adapted from cislunar logistics frameworks, to optimize resource generation, waste management, and disposal strategies. Future work will extend this analysis to include additional \ac{ISRU} methods, such as \acf{MCR} and \acf{HR}, which introduce unique challenges related to solid waste and emissions.

The paper is structured as follows: the Sec. \ref{sec:Methodology} details the framework for analyzing waste disposal logistics and \ac{ISRU} processes. Section \ref{sec:Results and Analysis} focuses on the comparison between Earth-dependent and ISRU-enabled mission scenarios, examining their economic and environmental implications. Finally, Section \ref{sec:Conclusion} concludes the significance of sustainable practices in lunar exploration and provides directions for future research.

\section{Methodology} \label{sec:Methodology}

\subsection{Objective Function and Constraints}

The methodology for this study revolves around a multi-commodity flow model that integrates the key aspects of resource generation, transportation logistics and waste disposal on the lunar surface. Particularly, the space mission planning formulation is represented as a multi-commodity flow logistics model, which is captured as a directed network graph \( G \). Consider a set of arcs \( \mathcal{A}= (\mathcal{N}, \mathcal{T}) \), which includes a set of nodes \( \mathcal{N} \) (indexed by \( i, j \)), and a set of time steps \( \mathcal{T} \) (indexed by \( t \)). 


The primary objective of the multi-commodity flow model is to minimize the total mission cost, $\mathcal{J}$, which includes transportation, resource, and waste generation costs. The objective function is defined as:

\begin{equation}
\textbf{Min} \quad
\mathcal{J} = \sum_{(i,j,t)\in\mathcal{A}} \bm{c}_{ij}^\top \bm{x}_{ijt}
\label{eq:objective_function}
\end{equation}

where:

\begin{itemize}[label={}]
    \item \( c_{ij} \): Cost coefficient for transporting commodities \( x_{ij} \).
    \item \( x_{ijt} \): Commodity flow variable representing mass transported from node \( i \) to \( j \) at time \( t \).
\end{itemize}

\subsubsection{Mass Balance Constraint}

Equation \ref{eq:mass_balance} represents the mass balance constraint, ensuring that the total inflow of commodities equals the sum of outflows and the demand at the node. Specifically, the second term \( Q_{ij} \) in the equation accounts for commodity transformations that occur either during spaceflights or after the deployment of infrastructure. These transformations include processes such as propellant consumption, the production of resources through \ac{ISRU} systems and the generation of byproducts or waste. The inflow and outflow calculations incorporate the interaction between various commodities, such as water, oxygen, hydrogen, and other byproducts. 

\begin{equation}
\sum_{(j):(i,j,t)\in\mathcal{A}} \bm{x}_{ijt} - \sum_{(j):(j,i,t)\in\mathcal{A}} Q_{ji}\bm{x}_{ji}(t-\Delta t_{ji}) = \bm{d}_{it} \quad \forall i \in \mathcal{N}, \forall t \in \mathcal{T}
\label{eq:mass_balance}
\end{equation}

where:
\begin{itemize}[label={}]
    \item \( Q_{ij} \): Transformation matrix for resource conversion or disposal processes during transport.
    \item \( d_{it} \): Demand (positive for resources, negative for waste) at node \( i \) at time \( t \).
    \item \( \Delta t_{ij} \): Time of flight or transportation duration between nodes \( i \) and \( j \).
\end{itemize}

\subsubsection{Concurrency Constraint}
This constraint limits the commodity flow based on various capacities, including the spacecraft's payload capacity, propellant capacity, \ac{ISRU} storage capacity, and waste disposal capacity. Equation \ref{eq:concurrency} ensures that the flow of commodities does not exceed the operational limits of the spacecraft, \ac{ISRU} systems, and disposal mechanisms. This constraint also guarantees that commodity inflows and outflows, including waste disposal flows, remain non-negative. Specifically, it ensures that the commodity flows are bounded by the available capacities and that no negative flows are allowed.

\begin{equation}
H_{ij}\bm{x}_{ijt} \leq \bm{0}_{l \times 1} \quad \forall (i,j,t) \in \mathcal{A}
\label{eq:concurrency}
\end{equation}

where:
\begin{itemize}[label={}]
    \item \( H_{ij} \): Capacity or operational matrix defining constraints on spacecraft, ISRU systems, or waste disposal operations.
    \item \( l \): Total number of concurrency constraint types.
\end{itemize}

\subsubsection{Time Disposal Window Constraint}
Equation \ref{eq:time_window} guarantees that commodity flows, including waste disposal flows, are only permitted when the corresponding time window is open for each arc between nodes. The time windows represent the intervals during which operations are allowed, and the constraint ensures that flows are zero outside of these valid time periods.

\begin{equation}
\begin{cases}
x_{ijt} = E_{\text{disposal}} & \text{if } t \in W_{ij} \\
x_{ijt} \geq 0, & \text{otherwise}
\end{cases}
\quad \forall k \in \mathcal{K} \quad \forall (i,j) \in \mathcal{A} \quad \forall t \in \mathcal{T}
\label{eq:time_window}
\end{equation}

where:
\begin{itemize}[label={}]
    \item \( W_{ij} \): Time Windows, which represent the valid intervals during which operations are allowed between the nodes \( i \) and \( j \).
\end{itemize}
The flow composition is defined as:

\begin{equation*}
\bm{x}_{ijt} = 
\begin{bmatrix}
\bm{x}_{\mathcal{C}} \\
\bm{x}_{D}
\end{bmatrix}_{ijt} \quad \bm{x}_{\mathcal{C}} \in \mathbb{R}^{|\mathcal{C}_c| \times 1}_{\geq 0}, \quad \bm{x}_{D} \in \mathbb{Z}^{|\mathcal{C}_D| \times 1}_{\geq 0}, \quad \forall (i,j,t) \in \mathcal{A}
\label{eq:flow_composition}
\end{equation*}

\subsection{Lunar Resources and Extraction} \label{Lunar Resources and Extraction.}

Lunar regolith is a layer of fragmented material produced by meteoritic impacts, solar wind, and extreme temperature cycling. Its composition is primarily silicates (plagioclase, pyroxene, and olivine) and oxide minerals such as ilmenite, which together constitute a rich source of oxygen bound in metal oxides\cite{schreiner2015} \cite{lu2011}. The regolith contains approximately 44\% oxygen by weight \cite{badescu2012, lu2011}, present in oxides of silicon, aluminum, iron and titanium. While its composition varies across the lunar surface, highlands are rich in plagioclase, whereas mare regions feature basaltic materials with higher iron and titanium content\cite{lu2011}. 

Also, water ice is now broadly acknowledged to exist on the lunar surface in specific locations, particularly in the Moon's permanently shadowed craters, which also house significant quantities of oxygen. Numerous studies \cite{anand2010, anand2012, basilevsky2012, crawford2015, hurley2016, li2018} have confirmed water ice presence, with concentrations estimated to range from 5\% to 30\% by weight in some areas. \cite{liu2023}\cite{schluter2020} This understanding was further reinforced in 2009 by the \ac{LCROSS} mission, which detected water ice during impact experiments in these cold, shadowed regions. The recovery of these water resources can support human activities through direct consumption, electrolysis for oxygen and hydrogen production, and other life-supporting processes\cite{fisher2010}.

Lunar \acf{ISRU} technologies enable the extraction and processing of local resources to support sustainable lunar exploration. These technologies reduce the dependency on Earth-based supplies by harnessing the Moon’s regolith and water ice for critical consumables such as oxygen, water, and hydrogen. The primary \ac{ISRU} methods for resource extraction are \acf{MRE}, \acf{SWE}, \acf{DWE}, \acf{MCR}, and \acf{HR}. Each method has unique chemical inputs and reactions, resulting in distinct byproducts that must be managed to prevent environmental contamination. This study focuses on analyzing the waste generation, disposal, and recycling potential of three of these methods: \ac{MRE}, \ac{SWE}, and \ac{DWE}. 

\begin{table}[hbt!]
\small 
\renewcommand{\arraystretch}{0.8} 
\centering
\caption{Lunar regolith processing methods.}
\begin{threeparttable}
\begin{tabular}{@{}p{1.7cm}p{2.5cm}p{2.5cm}p{2.3cm}p{2.8cm}p{2.2cm}@{}}
\toprule
\textbf{Process} & \textbf{Molten Regolith \newline Electrolysis \newline (MRE)} & \textbf{Methane \newline Carbothermal \newline Reduction (MCR)} & \textbf{Hydrogen \newline Reduction \newline (HR)} & \textbf{Soil/Water \newline Extraction \newline (SWE)} & \textbf{Direct Water \newline Electrolysis \newline (DWE)} \\ \midrule

\textbf{Applicable} & Lunar Mare \& \newline Highlands Regolith & Lunar Mare & Lunar Mare & Lunar Mare \& \newline Highlands Regolith & Water extracted from SWE \\

\textbf{Main Input\newline Materials} & Silicates (SiO\textsubscript{2}) \newline Oxides (FeO, TiO\textsubscript{2}) & Ilmenite (FeTiO\textsubscript{3})  Silicates (MgSiO\textsubscript{3}) & Ilmenite (FeTiO\textsubscript{3}) & Silicates (SiO\textsubscript{2})  & Water (H\textsubscript{2}O) \\

\textbf{Chemical \newline Reaction \cite{mueller2010}} & \footnotesize SiO\textsubscript{2} → Si + O\textsubscript{2} \newline 2FeTiO\textsubscript{3} → 2Fe +\newline  2TiO\textsubscript{2} + O\textsubscript{2} & \footnotesize FeTiO\textsubscript{3} + MgSiO\textsubscript{3} + CH\textsubscript{4} → CO + 2H\textsubscript{2} + Si + MgO & \footnotesize 2FeTiO\textsubscript{3} + 2H\textsubscript{2} → 2Fe + 2TiO\textsubscript{2} + 2H\textsubscript{2}O &\footnotesize Regolith → H\textsubscript{2}O + Dehydrated Soil & \footnotesize 2H\textsubscript{2}O→2H\textsubscript{2}+ O\textsubscript{2} \\

\textbf{Primary \newline Product} & Oxygen (O\textsubscript{2})  & Hydrogen(H\textsubscript{2}) & Water (H\textsubscript{2}O) & Water (H\textsubscript{2}O) & Oxygen (O\textsubscript{2}) \newline Hydrogen (H\textsubscript{2}) \\

\textbf{Byproducts \newline \& Waste} & Slag \newline Oxygen (O\textsubscript{2})\newline Silicon (Si)\newline Iron (Fe) \newline Titanium Dioxide (TiO\textsubscript{2}) 

& Carbon Monoxide (CO)\newline Hydrogen (H\textsubscript{2}) \newline Silicon (Si) \newline Magnesium Oxide (MgO) & Iron (Fe) \newline Titanium Dioxide (TiO\textsubscript{2})\newline Water (H\textsubscript{2}O) 

&  Dehydrated Soil \newline
Hydrogen sulfide \newline (H\textsubscript{2}S) \newline Ammonia (NH\textsubscript{3})\newline 
Sulfur dioxide (SO\textsubscript{2})
\newline ethylene (C\textsubscript{2}H\textsubscript{4})\newline Carbon dioxide (CO\textsubscript{2})\newline Methanol (CH\textsubscript{3}OH)\newline 
Methane (CH\textsubscript{4})\newline 
Hydroxyl (OH) & \makebox[1.8cm][c]{--} \\

\bottomrule
\end{tabular}
\end{threeparttable}
\end{table}

\subsubsection{\acf{MRE}}

\ac{MRE} is a technology that directly extracts oxygen and metals from lunar regolith by heating the material to 1600\degree and beyond. \cite{schreiner2015,sadoway2010,schreiner_parametric} This high-temperature process melts the regolith, allowing electrolysis to break down the metal oxides. Oxygen gas is released at the anode, while metals such as iron, silicon, aluminum, and titanium are deposited at the cathode. \cite{schreinerdevelopment} The absence of a need for chemical additives simplifies the operation, reducing the mass and complexity of the system.

The \ac{MRE} process is particularly appealing for lunar operations due to its high efficiency and resilience to regolith composition variations. Studies have shown that \ac{MRE} can extract up to 95\% of the oxygen contained in the regolith, achieving oxygen yields of approximately 28\% by mass. Additionally, the process has been validated for both highland and mare regolith, underscoring its robustness across different lunar terrains. \cite{schreiner_parametric}

Particularly, a 600 kg, 56.5 kW system can produce 10,000 kg of oxygen per year\cite{schreiner2015}, translating to 16.67 kg of oxygen per kilogram of system mass annually, highlighting the system's ability to meet high demand efficiently. \cite{schreiner2015}. Byproducts from the \ac{MRE} process include molten metals such as silicon, iron, and aluminum, as well as glassy slag. These materials can be repurposed for in-situ manufacturing, creating structural components or tools, thus minimizing waste. For example, silicon can be used for photovoltaic panels, while iron and aluminum can serve construction needs. The slag, which contains residual oxides, can be processed further or utilized in infrastructure development.\cite{landis2005, curreri2006,schreiner_parametric}

The amount of metal produced by an \ac{MRE} reactor depends on regolith type and operating temperature. For reactors producing 10,000 kg of oxygen per year, metal production increases with higher operating temperatures due to enhanced oxide reduction. In Mare regolith, the mass of leftover slag decreases by approximately 35\% as the temperature rises from 1850 K to 2300 K, with similar reductions observed for Highlands regolith at slightly lower peak temperatures of around 2000 K. \cite{schreiner_parametric} Reactors processing Highlands regolith can produce more silicon and aluminum, while Mare regolith yields higher quantities of iron and titanium.

At higher temperatures exceeding 2300 K,
the total metal production approaches the oxygen output, nearly matching the oxygen production level of 10,000 kg per year. \cite{schreiner_parametric} However, the metals tend to alloy at the bottom of the reactor, requiring post-processing for separation. Despite this additional step, the molten state of the alloy simplifies processing and minimizes energy expenditure. These metals, along with reduced slag production, enable the potential for \ac{MRE} reactors to support both resource extraction and material manufacturing.
Purified silicon extracted through \ac{MRE} can be refined using processes like plasma deposition and chemical vapor deposition for advanced applications such as photovoltaic cells, which are critical for lunar energy systems due to the Moon's high solar radiation density of approximately 1,370 W/m² compared to Earth's 950 W/m². \cite{liu2023} Additionally, metals like iron and aluminum can be utilized for constructing habitats, infrastructure, and scientific equipment, significantly enhancing the self-sufficiency of lunar operations. Slag, a byproduct of the process, can also be repurposed for stabilizing structures or creating radiation shields, \cite{sanders2011} minimizing waste and contributing to infrastructure development. 

\subsubsection{\acf{SWE} \& \acf{DWE}}

\ac{SWE} and \ac{DWE} are complementary methods focused on utilizing water sources from icy regolith or hydrated minerals. \ac{SWE} involves extracting water from the regolith, while \ac{DWE} subsequently electrolyzes this water to produce oxygen and hydrogen. \ac{SWE} produces non-volatile byproducts, such as slag processed soil with reduced water content—that requires careful disposal to prevent surface contamination. \ac{DWE}, by contrast, produces minimal solid waste, primarily generating gaseous byproducts like hydrogen, which can be stored or recycled within the \ac{ISRU} system. The integration of SWE and DWE promotes efficient use of lunar water resources, emphasizing the need for effective slag management to minimize environmental impact.

\ac{SWE} produces two primary byproducts: dehydrated soil and volatile emissions. Dehydrated soil accounts for 92.52\% of the regolith processed, making it the dominant byproduct. This material, significantly reduced in water content, may pose challenges for disposal or surface contamination if not managed effectively. Volatile emissions, on the other hand, constitute 1.88\% of the total regolith mass and are derived from approximately 33.49\% of the extracted water. These volatiles include compounds such as hydrogen sulfide (H\textsubscript{2}S), ammonia (NH\textsubscript{3}), sulfur dioxide (SO\textsubscript{2}), ethylene (C\textsubscript{2}H\textsubscript{4}), carbon dioxide (CO\textsubscript{2}), methanol (CH\textsubscript{3}OH), methane (CH\textsubscript{4}), and hydroxyl (OH). 

\ac{DWE} generates minimal solid waste. The main byproducts of \ac{DWE} are oxygen and hydrogen gases, which can be stored or directly utilized in mission operations. Oxygen constitutes 88.9\% of the water mass processed, while hydrogen accounts for the remaining 11.1\%. The efficiency of \ac{DWE} in separating these elements underscores its importance in providing breathable oxygen and energy-dense hydrogen for lunar operations.

\subsubsection{Excavator} 

The \ac{IPEx}  is a lightweight, 30-kg robotic system engineered to efficiently excavate and transport lunar regolith, processing up to 10 metric tons annually.\cite{schluter2020} Equipped with bucket drum technology, the excavator achieves an excavation rate of 333.33 kilograms, ensuring a continuous supply of raw materials for \ac{ISRU} processes. This advanced capability supports oxygen and water extraction, construction material generation, and sustainable operations on the Moon, significantly reducing reliance on Earth-launched resources and paving the way for long-term lunar infrastructure development.

\subsubsection{\acf{FSPS}}
The \ac{FSPS}, developed under NASA's Kilopower project, provides a reliable 10-kilowatt-class power source for critical lunar activities such as \ac{ISRU}, habitat operations, and scientific research \cite{nasa2024}. Its compact, nuclear fission-based design operates effectively in the Moon's harsh environment, including the prolonged lunar night. By reducing reliance on solar energy and large battery systems, the \ac{FSPS} serves as a scalable, long-term solution essential for sustainable lunar bases.

\begin{center}
\begin{table}[hbt]
\centering
\renewcommand{\arraystretch}{1.5}
\caption{Productivity for Each Component}
\begin{tabular}{p{2.9cm} p{1cm} p{4cm} p{1.8cm} >{\centering\arraybackslash}p{2.5cm} @{\,} >{\centering\arraybackslash}p{2.7cm}}
\toprule
\textbf{Reactor} & \textbf{Type} & \textbf{Product(P)/\newline Byproduct(Byp)} & \textbf{Mass \newline Fraction (\%)} & \textbf{Production Rate \newline (kg product/year\newline /kg plant)} & \textbf{Power \newline Consumption (kW/kg plant)} \\
\midrule

\parbox[t]{4cm}{\textbf{Molten Regolith \newline Electrolysis (MRE)}} 
& P & Oxygen (O$_2$), \( \alpha_{\text{MRE}}^{\text{O}_2} \) & 44\% \cite{badescu2012, lu2011} & 16.67 \cite{schreiner2015} & \parbox[t]{3cm}{\centering $-P^{\text{MRE}}$ \\ -0.0942 \cite{schreiner2015}}
 \\ \cmidrule(lr){3-5}

& \multirow{2}{*}{Byp} & Metals, \( \alpha_{\text{MRE}}^{\text{Metals}} \) &  40\% \cite{schreiner_parametric} & 15.16 & \\ 
\cmidrule(lr){3-5}

& & \hspace{0.5cm} Silicon (Si)  &  \hspace{0.05cm} 19.61\% \cite{schreiner_parametric}  & \hspace{0.1cm} 7.43  & \\
& & \hspace{0.5cm} Iron (Fe)     & \hspace{0.1cm} 9.80\%\cite{schreiner_parametric}  & \hspace{0.1cm} 3.71  & \\ 
& & \hspace{0.5cm} Aluminum (Al) & \hspace{0.1cm} 5.49\%\cite{schreiner_parametric}  & \hspace{0.1cm} 2.08  & \\ 
& & \hspace{0.5cm} Titanium (Ti) & \hspace{0.2cm}5.10\%\cite{schreiner_parametric}  & \hspace{0.1cm} 1.94  & \\ \cmidrule(lr){3-5}
& & Slag, \( \alpha_{\text{MRE}}^{\text{Slag}} \) & 16\% & 6.06  & \\ \cmidrule(lr){3-5}
&   & Total Dehydrated \newline Soil for MRE, \( -\beta_{\text{MRE}}^{\text{dsoil}} \) & -100\% & -37.89 & \\ 
\midrule

\parbox[t]{4cm}{\textbf{Soil/Water \newline Extraction (SWE)}} 
& P & Water (H$_2$O), \( \alpha_{\text{SWE}}^{\text{H}_2\text{O}} \) & 5.6\% \cite{colaprete2010}& 10.50 \cite{chen2020} & \parbox[t]{3cm}{\centering $-P^{\text{SWE}}$ \\ -0.0359 \cite{chen2020}} \\ \cmidrule(lr){3-5}

& \multirow{2}{*}{Byp} & Dehydrated Soil, \( \alpha_{\text{SWE}}^{\text{dsoil}} \) & 92.52\% & 173.48  & \\
& & Volatile Emissions, \( \alpha_{\text{SWE}}^{\text{emissions}} \) & 1.88\% \cite{colaprete2010} & 3.52 & \\ \cmidrule(lr){3-5}
& & Total Soil for SWE, \( -\beta_{\text{SWE}}^{\text{soil}} \) & -100\% & -187.5 & \\ 
\midrule

\parbox[t]{4cm}{\textbf{Direct Water \newline Electrolysis (DWE)}} 
& \multirow{3}{*}{P} & Oxygen (O$_2$), \( \alpha_{\text{DWE}}^{\text{O}_2} \) & 88.9\% & 31.12 \cite{chen2020} & \parbox[t]{3cm}{\centering $-P^{\text{DWE}}$ \\ -0.0700 \cite{chen2020}} \\ 
& & Hydrogen (H$_2$), \( \alpha_{\text{DWE}}^{\text{H}_2} \) & 11.1\% & 3.88 & \\ 
\cmidrule(lr){3-5}
&  & Total Water for DWE, \( -\beta_{\text{DWE}}^{\text{H}_2\text{O}} \) & -100\% & -35 & \\ 
\midrule

\parbox[t]{4cm}{\textbf{Excavator}} 
& P & Soil, \( \alpha_{\text{excavator}}^{\text{soil}} \)  & 100\% & 333.33 \cite{schuler2024ipex} &  \parbox[t]{3cm}{\centering \(  -P^{\text{excavator}} \) \\   -0.00113\cite{chen2020}}\\
\midrule

\parbox[t]{4cm}{\textbf{Fission Surface \newline Power System (FSPS)}} 
& P & Power Output & -- & -- &  \parbox[t]{3cm}{\centering \( P_{\text{FSPS}} \) \\ 0.00667 \cite{chen2020}}
\\ 
\bottomrule
\end{tabular}
\end{table}
\end{center}

\subsection{Lunar Surface Commodity Flow Network for Disposal Strategies} \label{Lunar Surface Commodity Flow Network for Disposal Strategies}

A key feature of this model is its ability to synthesize and quantify resource generation, transportation, and disposal over the full mission lifecycle. By tracking the flow of commodities into and out of the lunar surface system, the model enables the quantification of both economic and environmental impacts. This network-based lunar surface logistics methodology is designed to optimize the multi-commodity flow for post-mission disposal, focusing on evaluating the impacts of Artemis' key mission elements on the lunar surface. 

The multi-commodity flow logistics model is represented as a directed network graph, denoted as \( G \). The graph comprises a set of nodes, \( N \), and a set of arcs, \( A \), which capture the movement of mission-related commodities. In this model, nodes represent locations on Earth, \ac{LEO} and on the lunar surface, while arcs depict the possible transport and transformation pathways of mission elements over time. The proposed methodology tracks the commodity flows of both mission-critical resources and waste disposal elements during the life cycle of Artemis mission elements. 

This network is designed to facilitate the flow of commodities essential to lunar operations, including resources like oxygen (\(\text{O}_2\)), water (\(\text{H}_2\text{O}\)), hydrogen (\(\text{H}_2\)), as well as waste products like volatile emissions, processed slag and metals. Concretely, the commodity flow (\(x_{ijt}\))  represents the flow of a commodity \( c \in C \) from node \( i \) to node \( j \) at time step \( t \).

\[
X_{ijt} = \begin{bmatrix}
    x_{\text{soil}}, x_{\text{slag}}, x_{\text{dsoil}}, x_{\text{spares}}, x_{O_2}, x_{H_2}, x_{H_2O}, x_{\text{emissions}}, x_{\text{metals}}, x_{\text{excavator}}, x_{\text{SWE}}, x_{\text{DWE}}, x_{\text{MRE}}, x_{\text{power}}, x_{\text{FSPS}}
\end{bmatrix}^{T}_{ijt}
\]

 The transportation arcs in the network depict resource exchanges between different nodes, while the transformation arcs handle resource conversions and processes, such as \ac{ISRU}. The transformation matrix, $Q_{ij}$, allows for the computation of resource consumption and generation during lunar surface activities, similar to the use of propellant during spaceflights.

The lunar logistics model simulates the flow of multiple commodities through different mission stages, from resource extraction to waste disposal. The commodities considered include critical resources such as H\(_2\)O, O\(_2\), and H\(_2\), as well as infrastructure components like Excavators, Power Systems, and the Habitat. In addition, waste products such as processed slag, volatiles and emissions and other byproducts of ISRU activities are also included in the model.

The flow of the commodities is expressed through a system of equations, using a transformation matrix $Q_{ij}$, that tracks resource generation, conversion, and consumption. The following matrix provides an expanded overview of commodity interactions during lunar operations, incorporating not only ISRU systems but also waste disposal elements such as slag, metals, and gas emissions:

\begin{center}
\vspace{0.5cm}
\centerline{$
\renewcommand{\arraystretch}{1.1} 
\left[\begin{array}{@{}c@{}} 
x_{\text{soil}}       \\ [6pt]
x_{\text{slag}}       \\ [6pt] 
x_{\text{dsoil}}      \\ [6pt]
x_{\text{spares}}     \\ [6pt]
x_{\text{O}_2}        \\ [6pt]
x_{\text{H}_2}        \\ [6pt]
x_{\text{H}_2O}       \\ [6pt]
x_{\text{emissions}}  \\ [6pt]
x_{\text{metals}}     \\ [6pt]
x_{\text{excavator}}  \\ [6pt]
x_{\text{SWE}}        \\ [6pt]
x_{\text{DWE}}        \\ [6pt]
x_{\text{MRE}}        \\ [6pt]
x_{\text{power}}      \\ [6pt]
x_{\text{FSPS}}       \\ [6pt]
\end{array} \right]^{\text{inflow}}_{ijt} \hspace{-2mm} =
\left[\begin{array}{cccccccccc@{\hspace{0mm}}c@{\hspace{1mm}}c@{\hspace{0.5mm}}ccc}

1 & 0 & 0 & 0 & 0 & 0 & 0 & 0 & 0 & \alpha_{\text{excavator}}^{\text{soil}} & -\beta_{\text{SWE}}^{\text{soil}} & 0 & 0 & 0 & 0 \\[6pt]

0 & 1 & 0 & 0 & 0 & 0 & 0 & 0 & 0 & 0 & 0 & 0 & \alpha_{\text{MRE}}^{\text{slag}} & 0 & 0\\[6pt]

0 & 0 & 1 & 0 & 0 & 0 & 0 & 0 & 0 & 0 & \alpha_{\text{SWE}}^{\text{dsoil}} & 0 & -\beta_{\text{MRE}}^{\text{dsoil}} & 0 & 0\\[6pt]

0 & 0 & 0 & 1 & 0 & 0 & 0 & 0 & 0 & 0 & 0 & 0 & 0 & 0 & 0\\[6pt]

0 & 0 & 0 & 0 & 1 & 0 & 0 & 0 & 0 & 0 & 0 & \alpha_{\text{DWE}}^{\text{O}_2} & \alpha_{\text{MRE}}^{\text{O}_2} & 0 & 0\\[6pt]

0 & 0 & 0 & 0 & 0 & 1 & 0 & 0 & 0 & 0 & 0 & \alpha_{\text{DWE}}^{\text{H}_2} & 0 & 0  & 0\\[6pt]

0 & 0 & 0 & 0 & 0 & 0 & 1 & 0 & 0 & 0 & \alpha_{\text{SWE}}^{\text{H}_2O} & -\beta_{\text{DWE}}^{\text{H}_2O} & 0 & 0 & 0 \\[6pt]

0 & 0 & 0 & 0 & 0 & 0 & 0 & 1 & 0 & 0 & \alpha_{\text{SWE}}^{\text{emissions}} & 0 & 0 & 0 & 0 \\[6pt]

0 & 0 & 0 & 0 & 0 & 0 & 0 & 0 & 1 & 0 & 0 & 0 & \alpha_{\text{MRE}}^{\text{metals}} & 0 & 0\\[6pt]

0 & 0 & 0 & 0 & 0 & 0 & 0 & 0 & 0 & 1 & 0 & 0 & 0 & 0 & 0\\[6pt]

0 & 0 & 0 & 0 & 0 & 0 & 0 & 0 & 0 & 0 & 1 & 0 & 0 & 0 & 0\\[6pt]

0 & 0 & 0 & 0 & 0 & 0 & 0 & 0 & 0 & 0 & 0 & 1 & 0 & 0 & 0\\[6pt]

0 & 0 & 0 & 0 & 0 & 0 & 0 & 0 & 0 & 0 & 0 & 0 & 1 & 0 & 0\\[6pt]

0 & 0 & 0 & 0 & 0 & 0 & 0 & 0  & 0 & -P^{\text{excavator}} & -P^{\text{SWE}} & -P^{\text{DWE}} & -P^{\text{MRE}} & 1 & P_{\text{FSPS}}\\[6pt]

0 & 0 & 0 & 0 & 0 & 0 & 0 & 0 & 0 & 0 & 0 & 0 & 0 & 0 & 1\\[6pt]

\end{array} \right]_{ij}
\left[\begin{array}{@{}c@{}} 
x_{\text{soil}}       \\ [6pt]
x_{\text{slag}}       \\ [6pt] 
x_{\text{dsoil}}      \\ [6pt]
x_{\text{spares}}     \\ [6pt]
x_{\text{O}_2}        \\ [6pt]
x_{\text{H}_2}        \\ [6pt]
x_{\text{H}_2O}       \\ [6pt]
x_{\text{emissions}}  \\ [6pt]
x_{\text{metals}}     \\ [6pt]
x_{\text{excavator}}  \\ [6pt]
x_{\text{SWE}}        \\ [6pt]
x_{\text{DWE}}        \\ [6pt]
x_{\text{MRE}}        \\ [6pt]
x_{\text{power}}      \\ [6pt]
x_{\text{FSPS}}       \\ [6pt]
\end{array} \right]_{ijt} 
$}
\label{eq:tranformation_matrix} 
\end{center}

\section{Results and Analysis} \label{sec:Results and Analysis}

\subsection{Mission Scenario}

The mission planning framework in this study examines two logistical scenarios, focusing on their ability to meet resource demands for sustainable lunar operations: the Earth-Dependent Scenario and the ISRU-Enabled Scenario. These scenarios operate within a transportation network that includes Earth, \acf{LEO}, the \acf{LS}, and the \acf{HB} node. Both approaches aim to address annual demands for oxygen, hydrogen, and water while optimizing cost and environmental sustainability.

In the Earth-Dependent Scenario, all necessary resources are supplied from Earth. The annual demand at the Habitat Node consists of 10,000 kg of O\textsubscript{2} \cite{curreri2023} and 5,000 kg of H\textsubscript{2}O \cite{Lynch2023}, fully met through Earth-based logistics. These commodities are transported at a launch cost of \$5,000/kg. Liquid oxygen (LO\textsubscript{2}) costs \$0.15/kg, and liquid hydrogen (LH\textsubscript{2}) costs \$5.97/kg. Propellant consumption is calculated with a specific impulse (\(I_{\text{sp}}\)) of 420 seconds for LO\textsubscript{2}/LH\textsubscript{2} propulsion systems.

This scenario represents a straightforward logistical approach with 100\% reliance on Earth resources. However, the costs of transportation dominate the budget, driven by high launch expenses and the need for frequent resupply missions. Moreover, scalability is constrained by launch capacity, making it less viable for extended lunar operations.

\begin{figure}[hbt!]
\centering
\includegraphics[width=0.7\textwidth]{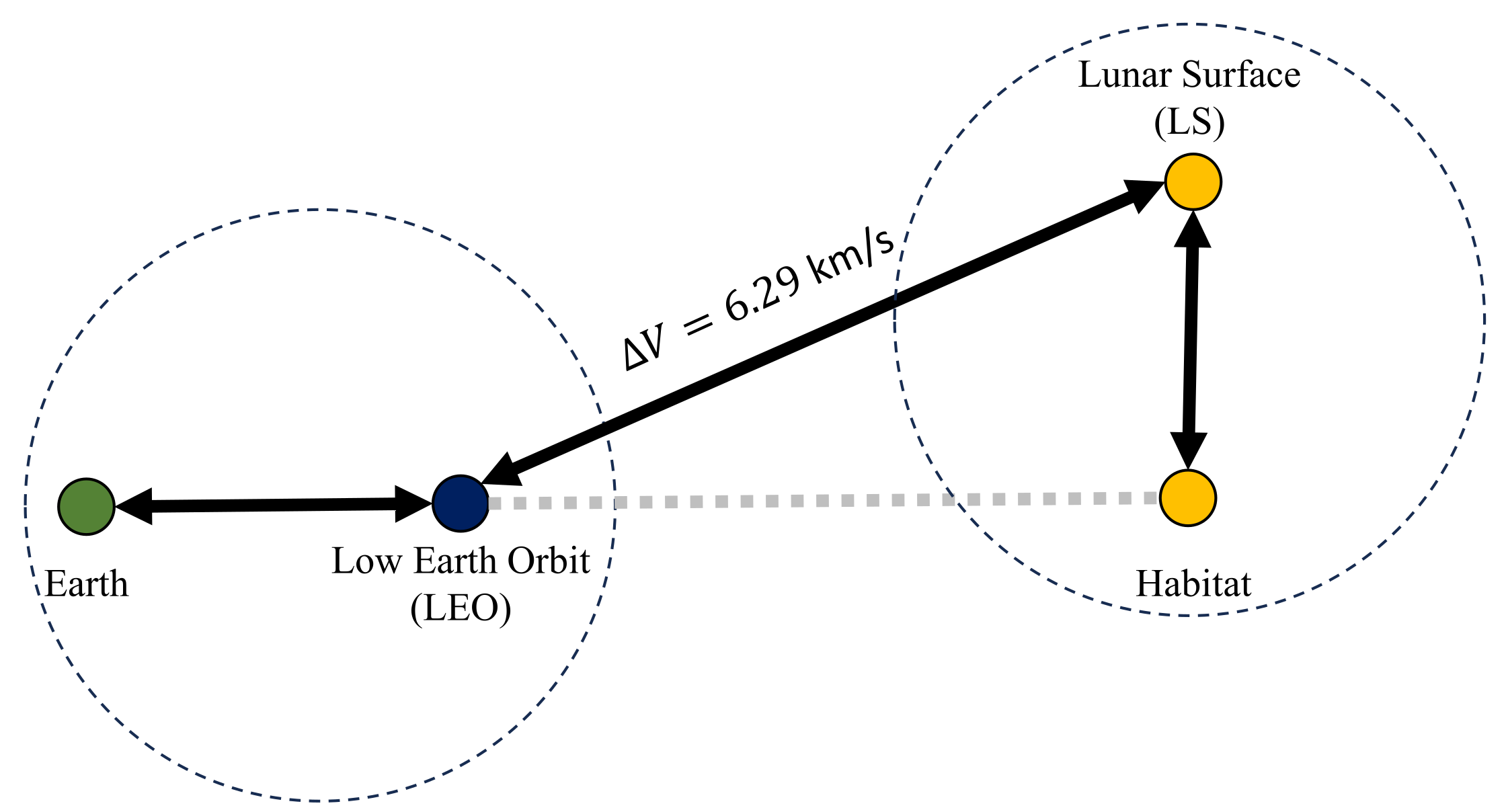}
\caption{Transportation Network Model.}
\label{fig:cislunar_network}
\end{figure}

Figure \ref{fig:cislunar_network} illustrates the multi-commodity flow structure across the lunar surface. In this framework, the lunar surface node is a synthetic interaction point between the different mission elements and the lunar environment. 

\begin{table}[hbt!]
\small 
\renewcommand{\arraystretch}{1.1} 
\centering
\caption{Assumptions of Mission Operation}
\label{tab: Assumptions of Mission Operation}
\begin{threeparttable}
\begin{tabular}{@{}p{4cm}p{2.5cm}@{}}
\toprule
\textbf{Parameter}                  & \textbf{Assumed Value} \\ \midrule
Propellant                          & LO\textsubscript{2}/LH\textsubscript{2}           \\ 
\(I_{\text{sp}}\)                   & 420 s       \\ 
Spacecraft propellant capacity      & 65,000 kg   \\ 
Spacecraft structure mass           & 6,000 kg    \\ 
ISRU maintenance                    & 5\% /year   \\ 
Rocket launch cost                  & \$5,000/kg   \\ 
Spacecraft manufacturing cost       & \$150 M      \\ 
Spacecraft Operation cost           & \$0.5 M      \\
LH\textsubscript{2} cost on Earth                   & \$5.97/kg            \\ 
LO\textsubscript{2} cost on Earth                   & \$0.15/kg            \\ 
Annual oxygen demand                & 10,000 kg \cite{curreri2023}\\ 
Annual water demand  (4 crew)       & 5,000 kg \cite{Lynch2023} \\ 
\bottomrule
\end{tabular}
\end{threeparttable}
\end{table}

\subsection{Mission Planning Results} \label{Mission Planning Results} 

The mission planning results, as summarized in Table \ref{tab:Comparative Analysis}, compare the cost and environmental impact of the ISRU mission with the Earth-dependent mission scenario. The mission assumes the parameters outlined in Table \ref{tab: Assumptions of Mission Operation}, including resource demands, technological capabilities, and operational constraints. These parameters provide the foundation for evaluating the comparative advantages of implementing \ac{ISRU} technologies.

The ISRU mission demonstrates a significant economic advantage, with a total cost of \$753M compared to \$2,481M for the Earth-dependent mission, making it approximately three times less expensive. This remarkable cost reduction is achieved through in-situ resource utilization, which enables the production of essential resources like oxygen and water directly on the lunar surface. In contrast, the Earth-dependent scenario requires launching all necessary materials from Earth, dramatically increasing mission costs due to higher launch mass and associated logistical challenges.

\begin{center}
\begin{table}[hbt!]
\centering
\renewcommand{\arraystretch}{1.2}
\caption{Comparative Analysis}
\label{tab:Comparative Analysis}
\begin{tabular}{p{3cm} p{3cm} p{3cm}}
\toprule
\textbf{Aspect}   & \textbf{ISRU Mission} & \textbf{Earth-Dependent} \\ \midrule
Total Cost                & \$753 M               & \$2,481M \\ 
Slag Produced             & 10,905.82 kg           & 0.0 kg     \\ 
Metals Produced           & 27,282.54 kg           & 0.0 kg     \\ 
Emissions                 & 5,028.57 kg            & 0.0 kg     \\ 

\bottomrule
\end{tabular}
\end{table}
\end{center}

\begin{figure}[hbt!]
\centering
\includegraphics[width=\textwidth]{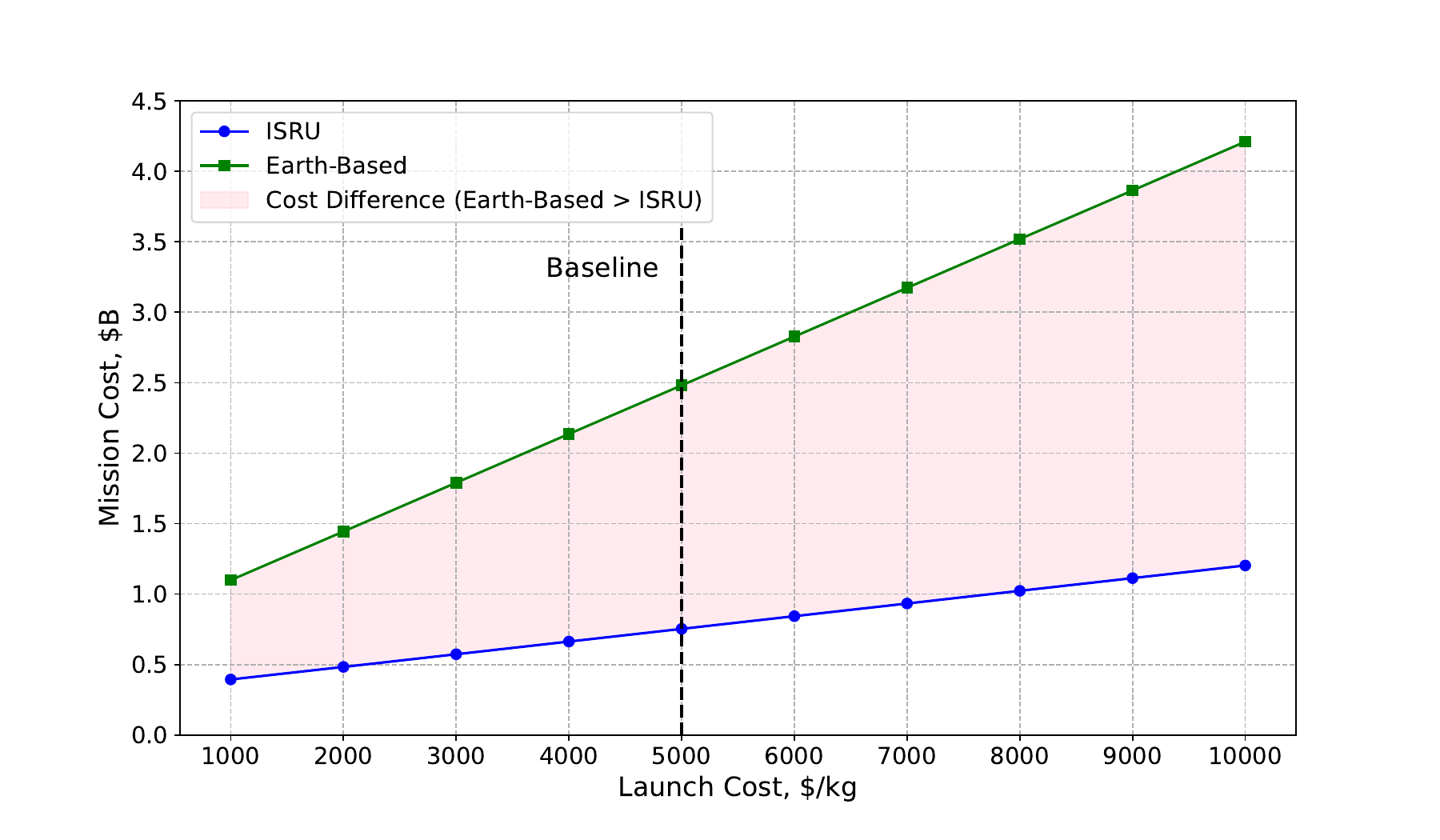}
\caption{Launch Cost vs Mission Cost.}
\label{fig:launch cost vs mission cost}
\end{figure}

Figure \ref{fig:launch cost vs mission cost} explores the relationship between launch cost and overall mission cost. It reveals a steady rise in mission cost for Earth-dependent missions as launch costs increase, reflecting their reliance on material transportation from Earth. Conversely, ISRU mission costs remain relatively stable, emphasizing the long-term cost-efficiency of in-situ resource utilization, particularly as launch costs escalate.

\begin{figure}[hbt!]
\centering
\includegraphics[width=\textwidth]{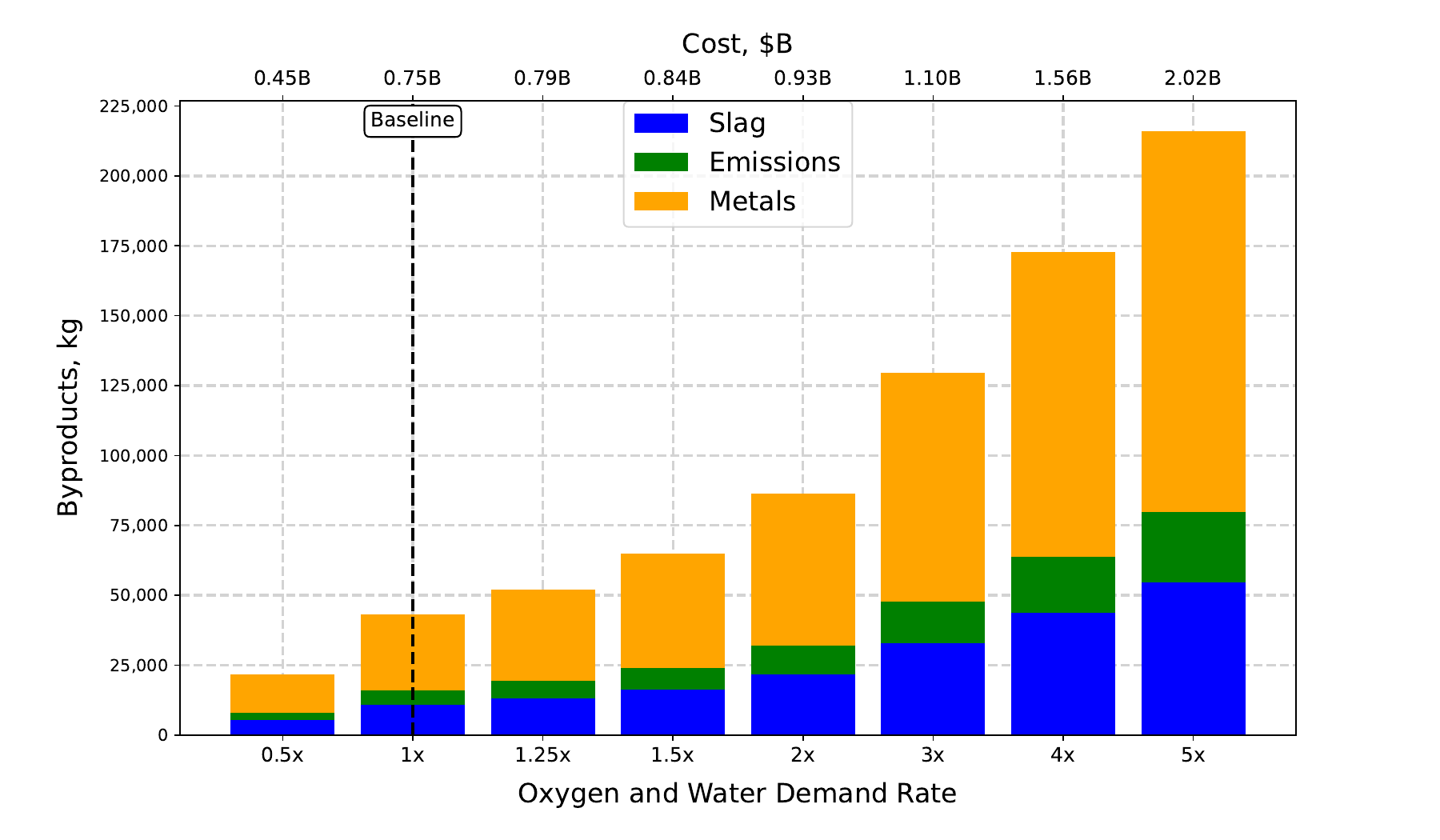}
\caption{ISRU Byproducts vs Productivity with Cost}
\label{fig:Byproducts vs Productivity with Cost}
\end{figure}

Figure \ref{fig:Byproducts vs Productivity with Cost} presents the effect of varying productivity rates (oxygen and water demand rates) on byproduct generation and mission costs. As productivity scales from 0.5x to 5x of the baseline, which is 10,000 kg oxygen and 5,000 kg water annual demand, the production of slag, metals, and emissions increases proportionally. This observation indicates a trade-off between higher resource generation rates and the associated increase in byproducts and environmental impact. Cost implications, shown in the figure, emphasize the economic feasibility of ISRU even at higher demand rates.

The comparative results underline the economic advantages of adopting ISRU technologies for lunar exploration. However, the environmental implications, particularly the emissions and slag byproducts, necessitate further innovation to mitigate adverse effects. These findings provide a foundation for optimizing \ac{ISRU} technologies to balance resource generation, cost-efficiency, and environmental impact.

\section{Conclusion} \label{sec:Conclusion}
This paper proposes a multicommodity flow-based approach to post-mission disposal and resource sustainability for lunar exploration. By integrating \ac{ISRU} systems, the proposed strategy demonstrates significant cost efficiency, achieving a threefold reduction in cost compared to Earth-dependent methods, with a total cost of \$753M versus \$2,481M. This approach leverages local resources, reducing the need for Earth-supplied materials while generating byproducts such as metals, volatile emissions and slag.

In space, waste retains potential economic value, particularly in extreme environments. Efficient ISRU practices not only support scalable exploration but also reduce the logistical burden of waste transportation and storage, paving the way for sustainable lunar bases. However, the byproducts of \ac{ISRU} must be carefully managed to mitigate environmental and contamination risks. Recycling, upcycling, and reuse of waste materials are imperative to minimize ecological impact and maximize resource utility.

\clearpage
\newpage
\bibliography{bibliography}

\begin{thebibliography}{30}
\newcommand{\enquote}[1]{``#1''}
\providecommand{\natexlab}[1]{#1}
\providecommand{\url}[1]{\texttt{#1}}
\providecommand{\urlprefix}{URL }
\expandafter\ifx\csname urlstyle\endcsname\relax
  \providecommand{\doi}[1]{\discretionary{}{}{}https://doi.org/#1}\else
  \providecommand{\doi}[1]{\discretionary{}{}{}\urlstyle{rm}\url{https://doi.org/#1}}\fi

\bibitem[{Zhang et~al.(2023)}]{zhang2023}
Zhang, P., et~al., \enquote{Overview of the Lunar In Situ Resource Utilization Techniques for Future Lunar Missions,} \emph{Space: Science \& Technology}, Vol. 2023, 2023, p.~37.
\newblock \doi{10.34133/space.0037}.

\bibitem[{NASA(2020)}]{nasa2020}
NASA, \enquote{NASA's Plan for Sustained Lunar Exploration and Development,} Tech. rep., National Aeronautics and Space Administration, 2020.
\newblock Retrieved from \url{https://www.nasa.gov/sites/default/files/atoms/files/a_sustained_lunar_presence_nspc_report4220final.pdf}.

\bibitem[{Anand(2010)}]{anand2010}
Anand, M., \enquote{Lunar Water: A Brief Review,} \emph{Earth, Moon, and Planets}, Vol. 107, No.~1, 2010, pp. 65--73.
\newblock \doi{10.1007/s11038-010-9377-9}.

\bibitem[{Badescu(2012)}]{badescu2012}
Badescu, V., \emph{Moon: Prospective Energy and Material Resources}, Springer, 2012.
\newblock Cited on pp. 22, 29.

\bibitem[{Hightower(1992)}]{hightower1992}
Hightower, T.~M., \enquote{Recycling and Source Reduction for Long Duration Space Habitation,} Tech. rep., NASA Johnson Space Center, July 1992.
\newblock \doi{10.4271/921121}.

\bibitem[{Sanders and Larson(2012)}]{sanders2012}
Sanders, G.~B., and Larson, W.~E., \enquote{Progress Made in Lunar In Situ Resource Utilization under NASA's Exploration Technology and Development Program,} Tech. rep., NASA Johnson Space Center, April 2012.
\newblock \doi{10.1061/9780784412190.050}.

\bibitem[{Basilevsky et~al.(2012)Basilevsky, Abdrakhimov, and Dorofeeva}]{basilevsky2012}
Basilevsky, A.~T., Abdrakhimov, A.~M., and Dorofeeva, V.~A., \enquote{Water and Other Volatiles on the Moon: A Review,} \emph{Solar System Research}, Vol.~46, No.~2, 2012, pp. 89--107.
\newblock \doi{10.1134/S0038094612010017}.

\bibitem[{Landis(2005)}]{landis2005}
Landis, G.~A., \enquote{Materials Refining for Solar Array Production on the Moon,} Tech. rep., NASA Technical Memorandum NASA/TM-2005-214014, 2005.
\newblock \urlprefix\url{https://ntrs.nasa.gov/citations/20050175853}.

\bibitem[{Curreri et~al.(2006)Curreri, Ethridge, Hudson, Miller, Grugel, Sen, and Sadoway}]{curreri2006}
Curreri, P., Ethridge, E., Hudson, S., Miller, T., Grugel, R., Sen, S., and Sadoway, D., \enquote{Process Demonstration for Lunar In Situ Resource Utilization: Molten Oxide Electrolysis,} Msfc independent research and development project, NASA Marshall Space Flight Center, 2006.
\newblock Project No. 5-81.

\bibitem[{Dallas et~al.(2021)Dallas, Raval, Saydam, and Dempster}]{dallas2021}
Dallas, J., Raval, S., Saydam, S., and Dempster, A., \enquote{An Environmental Impact Assessment Framework for Space Resource Extraction,} \emph{Space Policy}, 2021.
\newblock \doi{10.1016/j.spacepol.2021.101441}.

\bibitem[{Schreiner(2015)}]{schreiner2015}
Schreiner, S.~S., \enquote{Molten Regolith Electrolysis Reactor Modeling and Optimization of In-Situ Resource Utilization Systems,} Master of science in aerospace engineering thesis, Massachusetts Institute of Technology, Department of Aeronautical and Astronautical Engineering, June 2015.
\newblock Thesis Supervisor: Jeffrey A. Hoffman.

\bibitem[{Lu and Reddy(2011)}]{lu2011}
Lu, Y., and Reddy, R.~G., \enquote{Extraction of Metals and Oxygen from Lunar Soil,} \emph{High Temperature Materials and Processes}, Vol.~27, No.~4, 2011, pp. 223--232.
\newblock \doi{10.1515/HTMP.2008.27.4.223}.

\bibitem[{Anand et~al.(2012)Anand, Crawford, Balat‐Pichelin, Abanades, Van~Westrenen, Peraudeau, Jaumann, and Seboldt}]{anand2012}
Anand, M., Crawford, I., Balat‐Pichelin, M., Abanades, S., Van~Westrenen, W., Peraudeau, G., Jaumann, R., and Seboldt, W., \enquote{A Brief Review of Chemical and Mineralogical Resources on the Moon and Likely Initial In Situ Resource Utilization (ISRU) Applications,} \emph{Planetary and Space Science}, Vol.~74, No.~1, 2012, pp. 42--48.
\newblock \doi{10.1016/j.pss.2012.08.012}.

\bibitem[{Crawford(2015)}]{crawford2015}
Crawford, I.~A., \enquote{Lunar Resources: A Review,} \emph{Progress in Physical Geography}, Vol.~39, No.~2, 2015, pp. 137--167.
\newblock \doi{10.1177/0309133314567585}.

\bibitem[{Hurley et~al.(2016)Hurley, Colaprete, Elphic, Farrell, Hayne, Heldmann, Hibbits, Livengood, Lucey, and Klaus}]{hurley2016}
Hurley, D., Colaprete, A., Elphic, R., Farrell, W., Hayne, P., Heldmann, J., Hibbits, C., Livengood, T., Lucey, P., and Klaus, K., \enquote{Lunar Polar Volatiles: Assessment of Existing Observations for Exploration,} Tech. rep., NASA Goddard Space Flight Center, 2016.

\bibitem[{Li et~al.(2018)Li, Lucey, Milliken, Hayne, Fisher, Williams, Hurley, and Elphic}]{li2018}
Li, S., Lucey, P.~G., Milliken, R.~E., Hayne, P.~O., Fisher, E., Williams, J.-P., Hurley, D.~M., and Elphic, R.~C., \enquote{Direct Evidence of Surface Exposed Water Ice in the Lunar Polar Regions,} \emph{Proceedings of the National Academy of Sciences of the United States of America}, Vol. 115, No.~36, 2018, pp. 8907--8912.
\newblock \doi{10.1073/pnas.1802345115}.

\bibitem[{Liu et~al.(2023)Liu, Wang, Pang, Wang, Zhao, Lin, Wang, Shen, Liu, Song, Lai, Quan, and Yao}]{liu2023}
Liu, Y., Wang, C., Pang, Y., Wang, Q., Zhao, Z., Lin, T., Wang, Z., Shen, T., Liu, S., Song, J., Lai, X., Quan, X., and Yao, W., \enquote{Water Extraction from Icy Lunar Regolith by Drilling-Based Thermal Method in a Pilot-Scale Unit,} \emph{Acta Astronautica}, Vol. 202, 2023, pp. 386--399.
\newblock \doi{10.1016/j.actaastro.2022.11.002}.

\bibitem[{Schlüter and Cowley(2020)}]{schluter2020}
Schlüter, L., and Cowley, A., \enquote{Review of Techniques for In-Situ Oxygen Extraction on the Moon,} \emph{Planetary and Space Science}, Vol. 181, 2020, p. 104753.
\newblock \doi{10.1016/j.pss.2020.104753}.

\bibitem[{Fisher et~al.(2010)Fisher, Hogan, Pace, and Wignarajah}]{fisher2010}
Fisher, J.~W., Hogan, J.~A., Pace, G.~S., and Wignarajah, K., \enquote{Impact of Water Recovery from Wastes on the Lunar Surface Mission Water Balance,} \emph{40th International Conference on Environmental Systems Proceedings}, 2010.
\newblock \doi{10.2514/6.2010-6008}.

\bibitem[{Mueller et~al.(2010)Mueller, Townsend, Mantovani, and Metzger}]{mueller2010}
Mueller, R., Townsend, I., Mantovani, J., and Metzger, P., \enquote{Evolution of Regolith Feed Systems for Lunar ISRU O\textsubscript{2} Production Plants,} \emph{48th AIAA Aerospace Sciences Meeting Proceedings}, 2010.
\newblock \doi{10.2514/6.2010-1547}.

\bibitem[{Sadoway and Sibille(2010)}]{sadoway2010}
Sadoway, D.~R., and Sibille, L., \enquote{Direct Electrolysis of Molten Lunar Regolith for the Production of Oxygen and Metals on the Moon,} \emph{ECS Transactions}, 2010.
\newblock \doi{10.1149/1.3367929}.

\bibitem[{Schreiner et~al.(2016)Schreiner, Sibille, Dominguez, and Hoffman}]{schreiner_parametric}
Schreiner, S.~S., Sibille, L., Dominguez, J.~A., and Hoffman, J.~A., \enquote{A parametric sizing model for Molten Regolith Electrolysis reactors to produce oxygen on the Moon,} \emph{Advances in Space Research}, Vol.~57, No.~7, 2016, pp. 1585--1603.
\newblock \doi{https://doi.org/10.1016/j.asr.2016.01.006}.

\bibitem[{Schreiner et~al.(2015)Schreiner, Sibille, Dominguez, and Sirk}]{schreinerdevelopment}
Schreiner, S., Sibille, L., Dominguez, J., and Sirk, A., \enquote{Development of a Molten Regolith Electrolysis Reactor Model for Lunar In-Situ Resource Utilization,} \emph{8th Symposium on Space Resource Utilization Proceedings}, 2015.
\newblock \doi{10.2514/6.2015-1180}.

\bibitem[{Sanders(2011)}]{sanders2011}
Sanders, G.~B., \enquote{Comparison of Lunar and Mars In-Situ Resource Utilization for Future Robotic and Human Missions,} Tech. rep., NASA Johnson Space Center, January 2011.
\newblock \doi{10.2514/6.2011-120}.

\bibitem[{{NASA}(2024)}]{nasa2024}
{NASA}, \enquote{{Kilopower: NASA's Small Nuclear Power System for Space Exploration},} , 2024.
\newblock \urlprefix\url{https://www.nasa.gov/directorates/stmd/tech-demo-missions-program/kilopower-hmqzw/}.

\bibitem[{Colaprete et~al.(2010)Colaprete, Schultz, Heldmann, Wooden, Shirley, Ennico, and Sollitt}]{colaprete2010}
Colaprete, A., Schultz, P., Heldmann, J., Wooden, D., Shirley, M., Ennico, K., and Sollitt, L., \enquote{Detection of Water in the LCROSS Ejecta Plume,} \emph{Science}, Vol. 330, No. 6003, 2010, pp. 463--468.
\newblock \doi{10.1126/science.1186986}.

\bibitem[{Chen et~al.(2020)Chen, Sarton~du Jonchay, Hou, and Ho}]{chen2020}
Chen, H., Sarton~du Jonchay, T., Hou, L., and Ho, K., \enquote{Integrated In-Situ Resource Utilization System Design and Logistics for Mars Exploration,} \emph{Acta Astronautica}, Vol. 170, 2020, pp. 80--92.
\newblock \doi{10.1016/j.actaastro.2020.01.031}.

\bibitem[{Schuler et~al.(2024)Schuler, Nick, Leucht, Langton, and Smith}]{schuler2024ipex}
Schuler, J., Nick, A., Leucht, K., Langton, A., and Smith, D., \enquote{ISRU Pilot Excavator: Bucket Drum Scaling Experimental Results,} Tech. rep., NASA Kennedy Space Center, Kennedy Space Center, FL, USA, 2024.
\newblock \urlprefix\url{https://ntrs.nasa.gov/citations/20210025846}.

\bibitem[{Curreri(2023)}]{curreri2023}
Curreri, P., \enquote{Lunar South Pole Oxygen Pipeline,} , January 2023.
\newblock \urlprefix\url{https://www.nasa.gov/general/lunar-south-pole-oxygen-pipeline/}, published on NASA website.

\bibitem[{Lynch et~al.(2023)Lynch, Goodliff, Stromgren, Vega, and Ewert}]{Lynch2023}
Lynch, C.~S., Goodliff, K.~E., Stromgren, C., Vega, J., and Ewert, M.~K., \enquote{Logistics Rates and Assumptions for Future Human Spaceflight Missions Beyond LEO,} , 2023.
\newblock \urlprefix\url{https://ntrs.nasa.gov/citations/20230012635}, nASA Technical Report, Document ID: 20230012635.

\end{thebibliography}

\end{document}